# Modular Irregularity Strength of Triangular Book Graph

Meilin Imelda Tilukay

Department of Mathematics, Faculty of Mathematics and Natural Sciences, Pattimura University,
Jalan Ir. M. Putuhena, Kampus Poka – Unpatti, Ambon, Indonesia.
Email: meilinity@gmail.com

**Abstract:** This paper deals with the modular irregularity strength of a graph of $n$ vertices, a new graph invariant, modified from the irregularity strength, by changing the condition of the vertex-weight set associate to the well-known irregular labeling from $n$ distinct positive integer to $Z_n$-the group of integer modulo $n$. Investigating the triangular book graph $B_m^{(3)}$, we first find the irregularity strength of triangular book graph $s\left(B_m^{(3)}\right)$, as the lower bound for the modular irregularity strength, and then construct a modular irregular $s\left(B_m^{(3)}\right)$-labeling. The result shows that triangular book graphs admit a modular irregular labeling and its modular irregularity strength and irregularity strength are equal, except for a small case and the infinity property.

## 1. Introduction

For a simple graph $G$ of order $n \geq 2$, it is impossible to have $n$ distinct vertex degree. By adding multiple edges to $G$, each vertex can have distinct degree. It means that multigraph can have that property. A graph is irregular if its vertices have distinct degrees. Replacing multiple edges joining each pair of vertices by its number attached to an edge, Chartrand, Jacobson, Lehel, Oellermann, Ruiz, and Saba in [8] introduced the well-known labeling of $G$, called the *irregular assignment*, that is an edge $k$-labeling of the edge-set $f: E(G) \to \{1, 2, \cdots, k\}$ such that the vertex-weights are all distinct, where the weight of a vertex $u$ in $G$ is the sum of all labels of edges incident to $u$, wrote $w_f(u) = \sum_{uv \in E(G)} f(uv)$. Irregular assignment is also called a *vertex irregular edge $k$-labeling*. The minimum value $k$ for which $G$ has a vertex irregular edge $k$-labeling is called the *irregularity strength* of $G$, denoted by $s(G)$. If $G$ has no such labeling, $s(G) = \infty$. $s(G)$ is finite for only graph that contain no component of order at most 2.

The lower bound of $s(G)$ is given in [8] as follow.

$$s(G) \geq \max_{1 \leq i \leq \Delta} \left\{ \frac{n_i + i - 1}{i} \right\}, \tag{1}$$

where $n_i$ denotes the number of vertices of degree $i$ and $\Delta$ is the maximum degree of $G$. For $r$-regular graphs of order $n$, the lower bound [8] is $s(G) \geq \frac{n+r-1}{r}$.

In [9], Faudree and Lehel provided the upper bound for $r$-regular graphs of order $n$, $r \geq 2$, as $s(G) \leq \left\lceil \frac{n}{2} \right\rceil + 9$, and conjectured that $s(G) \leq \frac{n}{d} + c$ for any graph. For even $r$, Faudree, Schelp, Jacobson, and Lehel in [10] proved that $s(G) \leq \left\lceil \frac{n}{2} \right\rceil + 2$. In [14], Nierhoff gave a tight bound for general graph of order $n$ with no component of order at most 2, as $s(G) \leq n - 1$. Kalkowski, Karonski, and Pfender [12] also improved the bound of $s(G)$. The exact values of the irregularity strength of graphs are known only for few family of graphs, such as paths and complete graphs [8], cycles and Turan graphs [10], circulant graphs [1], trees [7], corona product of path and complete graphs, and corona product of cycle and complete graphs for small cases [13], fan graphs [3], wheel graphs [6].

Difficulties of finding the irregularity strength for any graph or even for graphs with simple structure have brought out many modifications of such labeling that one can find in [2], [4], [15], [16], [17], and [11]. The recent one is a modular irregular labeling of a graph introduced by Bača, Muthugurupackiam, Kathiresan, and Ramya in [5], which is obtained by modifying the condition of the vertex-weight set associate to the irregular labeling from $n$ distinct positive integer to $Z_n$-the group of integer modulo $n$.

Let $G = (V, E)$ be a graph of order $n$ with no component of order at most 2. An edge $k$-labeling $f: E(G) \to \{1, 2, \cdots, k\}$ is called a *modular irregular $k$-labeling* of $G$ if there exists a bijective weight function $w_f: V(G) \to Z_n$ defined by

$$w_f(x) = \sum f(xy)$$

called *modular weight* of the vertex $x$, where $Z_n$ is the group of integers modulo $n$ and the sum is over all vertices $y$ adjacent to $x$. The minimum value $k$ for which $G$ admits a modular irregular $k$-labeling is called the *modular irregularity strength* of $G$ and denoted by $ms(G)$. If a graph $G$ admits no modular irregular $k$-labeling, then $ms(G) = \infty$.

The lower bound of the modular irregularity strength of a graph is given in [5] as follow.
$$ms(G) \leq s(G). \tag{2}$$
And for any graph of order $n$, the infinity condition is given in Theorem A.

**Theorem A** ([5]). If $G$ is a graph of order $n$, $n \equiv 2 \pmod{4}$, then $G$ has no modular irregular labeling i.e., $ms(G) = \infty$.

A condition for an irregular assignment of a graph $G$ is also its modular irregular labeling is given in Lemma B.

**Lemma B** ([5]). Let $G$ be a graph with no component of order $\leq 2$, and let $s(G) = k$. If there exists an irregular assignment of $G$ with edge labels of at most $k$, where the weights of vertices constitute a set of constitute integer, then
$$s(G) = ms(G) = k.$$

They [5] also determined the exact values of the modular irregularity strength of five families of graphs, such as paths, stars, triangular graphs, cycles, and gear graphs. Later in [3], Bača, Kimáková, Lasscsáková, and Semaničová-Feňovčiková determined the modular irregularity strength of fan graphs, and in [6], Bača, Imran, and Semaničová-Feňovčiková determined the modular irregularity strength of wheel graphs. They ([3] and [6]) proposed the following problem.

**Problem 1** ([6]). Is there another family of graphs, besides wheels and fan graphs, for which the irregularity strength and the modular irregularity strength are the same?

The triangular book graph $B_n^{(3)}$, $n \geq 1$, is a planar undirected graph of order $n + 2$ and size $2n + 1$ constructed by $n$ cycles of order 3 sharing a common edge.

In the next section, we discuss the irregularity strength and the modular irregularity strength of triangular book graphs, in order to provide small answer to the problem.

## 2. Main Results

The first result gives the exact value of $s\left(B_n^{(3)}\right)$, $n \geq 1$.

### 2.1. The Irregularity Strength of Triangular Book Graphs

*Theorem 1.* Let $B_n^{(3)}$, $n \geq 1$, be a triangular book graph of order $n + 2$ and size $2n + 1$. Then

$$s\left(B_n^{(3)}\right) = \begin{cases} 3, & \text{for } n = 1 \\ \left\lceil \dfrac{n+1}{2} \right\rceil, & \text{for } n \geq 2 \end{cases}.$$

*Proof.* Let $B_n^{(3)}$, $n \geq 1$, be a triangular book graph with the vertex set $V\left(B_n^{(3)}\right) = \{a, b, c_i | 1 \leq i \leq n\}$ and the edge set $E\left(B_n^{(3)}\right) = \{ab, ac_i, bc_i | 1 \leq i \leq n\}$. We divide the proof into 2 cases as follow.

*Case 1.* For $n = 1$. It is clear that $B_1^{(3)}$ isomorphic to a cycle $C_3$, then $B_1^{(3)}$ admits a vertex irregular 3-labeling with edge labels $1, 2, 3$, and the induced vertex weights $3, 4, 5$, and $s\left(B_1^{(3)}\right) = 3$.

*Case 2.* For $n \geq 2$. By equation (1), we have that since $B_n^{(3)}$ is a bicenter graph with $\delta = 2$, then $s\left(B_n^{(3)}\right) \geq \left\lceil \frac{n+1}{2} \right\rceil$. The sufficient condition to complete the proof is by constructing a vertex irregular edge $\left\lceil \frac{n+1}{2} \right\rceil$-labeling. Define a vertex irregular edge $\left\lceil \frac{n+1}{2} \right\rceil$-labeling $f: E\left(B_n^{(3)}\right) \to \left\{1, 2, \cdots, \left\lceil \frac{n+1}{2} \right\rceil\right\}$ as follow.

$$f(ab) = \begin{cases} 2, & \text{for } n = 2; \\ 1, & \text{for } n \geq 3; \end{cases}$$

$$f(ac_i) = \begin{cases} \dfrac{i+1}{2}, & \text{for odd } i; \\ \dfrac{i}{2}, & \text{for even } i; \end{cases}$$

$$f(bc_i) = \begin{cases} \dfrac{i+1}{2}, & \text{for odd } i; \\ \dfrac{i}{2} + 1, & \text{for even } i. \end{cases}$$

It is clearly to see that the maximum label is $\left\lceil \frac{n+1}{2} \right\rceil$, hence, $f$ is an edge $\left\lceil \frac{n+1}{2} \right\rceil$-labeling. Next, we evaluate the vertex-weights as follow.

For $n = 2$, we have $w_f(c_1) = 2$, $w_f(c_2) = 3$, $w_f(a) = 4$, and $w_f(b) = 5$.

For odd $n \geq 3$, we have

$w_f(c_i) = i + 1$, $1 \leq i \leq n$;

$w_f(a) = \frac{1}{4}(n^2 + 2n + 5)$;

$w_f(b) = \frac{1}{4}(n^2 + 4n + 3)$;

For even $n \geq 4$, we have

$w_f(c_i) = i + 1$, $1 \leq i \leq n$;

$w_f(a) = \frac{1}{4}(n^2 + 2n + 4)$;

$w_f(b) = \frac{1}{4}(n^2 + 4n + 4)$.

The labeling $f$ is optimal and the vertex weights are all distinct, with $w_f(c_i) < w_f(a) < w_f(b)$, hence, $f$ is a vertex irregular edge $\left\lceil \frac{n+1}{2} \right\rceil$-labeling. Then, it can be concluded that $B_n^{(3)}$ admits a vertex irregular $\left\lceil \frac{n+1}{2} \right\rceil$-labeling and the irregularity strength $s\left(B_n^{(3)}\right) = \left\lceil \frac{n+1}{2} \right\rceil$. ∎

## 2.2. The Modular Irregularity Strength of Triangular Book Graphs

*Theorem 2.* Let $B_n^{(3)}$, $n \geq 1$, be a triangular book graph of order $n+2$ and size $2n+1$. Then

$$ms\left(B_n^{(3)}\right) = \begin{cases} 3, & \text{for } n = 1; \\ 4, & \text{for } n = 5; \\ \infty, & \text{for } n \equiv 0 \pmod{4}; \\ \left\lceil \dfrac{n+1}{2} \right\rceil, & \text{otherwise} \end{cases}$$

*Proof.* Let $B_n^{(3)}$, $n \geq 1$, be a triangular book graph with the vertex set $V\left(B_n^{(3)}\right) = \{a, b, c_i | 1 \leq i \leq n\}$ and the edge set $E\left(B_n^{(3)}\right) = \{ab, ac_i, bc_i | 1 \leq i \leq n\}$. We divide the proof into 4 cases as follow.

*Case 1.* For $n = 1$. It follows from Theorem 1 and Lemma A that $ms\left(B_1^{(3)}\right) = 3$.

*Case 2.* For $n = 5$. By Theorem 1 and equation (2), we have $ms\left(B_5^{(3)}\right) \geq 3$. Suppose that $B_5^{(3)}$ admits a modular irregular 3-labeling $f$ and $ms\left(B_5^{(3)}\right) = 3$. Since the degree of $c_i$, $1 \leq i \leq 5$, is 2 then the vertex weight under labeling $f$ is at least 2 and at most 6. Then the modular weight 0 and 1, must be realizable on both centers $a$ and $b$. Moreover, when we set the vertex weights $2, 3, \cdots, 6$, we obtained that the minimum weight of vertex $a$ and $b$ is $1 + 1 + 1 + 1 + 2 + 3 = 9$, and the maximum one is $3 + 1 + 2 + 3 + 3 + 3 = 15$. Since, $14 \equiv 0 \pmod{7}$ and $15 \equiv 1 \pmod{7}$, then we need to have the weight of vertices $a$ and $b$ equal to 14 and 15, respectively. Assume that $wt(b) = f(ab) + f(bc_1) + f(bc_2) + \cdots + f(bc_5) = 15$, then the only solution for the weight of vertex $a$ is $wt(a) = f(ab) + f(ac_1) + f(ac_2) + \cdots + f(ac_5) = 3 + 1 + 1 + 1 + 2 + 3 = 11 \equiv 4 \pmod{7}$ equals to one of modular weight we have set, which is a contradiction. The proof is similar for taking $wt(a) = 14$. Thus, $ms\left(B_5^{(3)}\right) \geq 4$. Let the edge labels listed as $f(ab) = 1, f(ac_1) = 1, f(ac_2) = 1, f(ac_3) = 1, f(ac_4) = 2, f(ac_5) = 2, f(bc_1) = 1, f(bc_2) = 2, f(bc_3) = 3, f(bc_4) = 3, f(bc_5) = 4$ be the construction of a modular irregular 4-labeling of $B_5^{(3)}$ such that the modular weights obtained are $w_f(c_i) = i + 1$, $1 \leq i \leq 5$, $w_f(a) = 8 \equiv 1 \pmod 7$, and $w_f(b) = 14 \equiv 0 \pmod 7$.

*Case 3.* For $n \equiv 0 \pmod 4$, it follows from Theorem A that $ms\left(B_n^{(3)}\right) = \infty$.

*Case 4.* For $n \neq 1, 5$ and $n \not\equiv 0 \pmod 4$, by Theorem 1 and equation (2), we have $ms\left(B_n^{(3)}\right) \geq \left\lceil \frac{n+1}{2} \right\rceil$. Next, we construct a vertex irregular edge $\left\lceil \frac{n+1}{2} \right\rceil$-labeling. Define a vertex irregular edge $\left\lceil \frac{n+1}{2} \right\rceil$-labeling $f: E\left(B_n^{(3)}\right) \to \left\{1, 2, \cdots, \left\lceil \frac{n+1}{2} \right\rceil\right\}$ as follow.

For $n \equiv 1 \pmod 8$,
$f(ab) = 1$;

$$f(ac_i) = \begin{cases} 1, & \text{for } 1 \leq i \leq \frac{n-1}{2}; \\ \frac{n-1}{8} + 2, & \text{for } i = \frac{n+1}{2}; \\ \frac{2i - n + 1}{2}, & \text{for } \frac{n+3}{2} \leq i \leq n; \end{cases}$$

$$f(bc_i) = \begin{cases} i, & \text{for } 1 \leq i \leq \frac{n-1}{2}; \\ \frac{3n-3}{8}, & \text{for } i = \frac{n+1}{2}; \\ \frac{n+1}{2}, & \text{for } \frac{n+3}{2} \leq i \leq n. \end{cases}$$

For $n \equiv 5 \pmod 8$,
$f(ab) = 1$;

$$f(ac_i) = \begin{cases} 1, & \text{for } 1 \le i \le \frac{n-1}{2}; \\ \frac{n+1}{2}, & \text{for } i = \frac{n+1}{2}; \\ \frac{n+35}{8}, & \text{for } i = \frac{n+3}{2}; \\ \frac{2i-n+1}{2}, & \text{for } \frac{n+5}{2} \le i \le n; \end{cases}$$

$$f(bc_i) = \begin{cases} i, & \text{for } 1 \le i \le \frac{n-1}{2}; \\ 1, & \text{for } i = \frac{n+1}{2}; \\ \frac{3n-15}{8}, & \text{for } \frac{n+3}{2}; \\ \frac{2i-n+1}{2}, & \text{for } \frac{n+5}{2} \le i \le n. \end{cases}$$

For $n \equiv 2 \pmod 4$, $n \equiv 3 \pmod 4$, and $1 \le i \le n$.

$$f(ab) = \begin{cases} \frac{n+6}{4}, & \text{for } n \equiv 2 \pmod 4; \\ 1, & \text{for } n \equiv 3 \pmod 4; \end{cases}$$

$$f(ac_i) = \begin{cases} \frac{i+1}{2}, & \text{for odd } i; \\ \frac{i}{2} + 1, & \text{for } i \equiv 0 \pmod 4; \\ \frac{i}{2}, & \text{for } i \equiv 2 \pmod 4; \end{cases}$$

$$f(bc_i) = \begin{cases} \frac{i+1}{2}, & \text{for odd } i; \\ \frac{i}{2}, & \text{for } i \equiv 0 \pmod 4; \\ \frac{i}{2} + 1, & \text{for } i \equiv 2 \pmod 4. \end{cases}$$

It can be checked that the maximum label given above is $\left\lceil \frac{n+1}{2} \right\rceil$, hence, $f$ is an edge $\left\lceil \frac{n+1}{2} \right\rceil$-labeling. Next, we evaluate the vertex-weights as follow.

For $n \equiv 1 \pmod 8$, we have $wt(c_i) = i + 1$, $1 \le i \le n$; $wt(a) = \frac{1}{8}(n+7)(n+2) \equiv 0 \pmod{(n+2)}$; and $wt(b) = \frac{3}{8}(n-1)(n+2) + 1 \equiv 1 \pmod{(n+2)}$.

For $n \equiv 5 \pmod 8$, we have $wt(c_i) = i + 1$, $1 \le i \le n$; $wt(a) = \frac{1}{8}(n+11)(n+2) \equiv 0 \pmod{(n+2)}$; and $wt(b) = \frac{1}{8}(3n-7)(n+2) + 1 \equiv 1 \pmod{(n+2)}$.

For $n \equiv 2 \pmod 4$, we have $wt(c_i) = i + 1$, $1 \le i \le n$; $wt(a) = \frac{1}{4}(n+2)^2 \equiv 0 \pmod{(n+2)}$; and $wt(b) = \frac{1}{4}(n+2)^2 + 1 \equiv 1 \pmod{(n+2)}$.

For $n \equiv 3 \pmod 4$, we have $wt(c_i) = i + 1$, $1 \le i \le n$; $wt(a) = \frac{1}{4}(n+1)(n+2) \equiv 0 \pmod{(n+2)}$; and $wt(b) = \frac{1}{4}(n+1)(n+2) + 1 \equiv 1 \pmod{(n+2)}$.

The labeling $f$ is optimal and forms the modular weights set $\{0, 1, 2, \cdots, n+1\}$. It means that there exists a bijective weight function $wt: V\left(B_n^{(3)}\right) \to Z_{n+2}$, such that $f$ satisfy a modular irregular $\left\lceil \frac{n+1}{2} \right\rceil$-labeling. ∎

## 3. Conclusion

By Theorem 1 and Theorem 2, we have determined the exact values of irregularity strength and modular irregularity strength of triangular book graphs. We conclude that for $n \ne 5$ and $n \not\equiv 0 \pmod 4$, $s\left(B_n^{(3)}\right) = ms\left(B_n^{(3)}\right)$.